\documentstyle[12pt]{article}
\input{amssym.def}

\advance\skip\footins by 0.3cm
\makeatletter
\def\@biblabel#1{#1.}

\makeatletter
\ifnum\@ptsize=0  \fi
\ifnum\@ptsize=2  \fi
\advance\textheight by 0.7cm
\advance\voffset by -0.3cm
\newcommand{\proof}{\noindent{\em Proof.}\, }
\newcommand\al{\alpha}
\newcommand\as[1]{\alpha_{#1}}
\newcommand\be[1]{\beta_{#1}}
\newcommand{\qed}{\hspace*{\fill}Q.E.D.\vskip12pt plus 1pt}
\newcommand\pn[1]{{\bf P}^{#1}}
\newcommand\ol[1]{\overline{#1}}

\newcommand\cQ{{\cal Q}}
\newcommand\Qs[1]{{\cal Q}_{#1}}
\newcommand\Os[1]{{\cal O}_{#1}}
\newcommand\cO{{\cal O}}

\newcommand\cH{{\cal H}}
\newcommand\cF{{\cal F}}
\newcommand\sub{\subset}
\newcommand\Fu{F^{i,j}}
\newcommand\Fd{F_{i,j}}
\newtheorem{theorem}{{\rm T\sc heorem}}[subsection]
 \newtheorem{lemma}[theorem]{{\rm L\sc emma}}
 
\newtheorem{proposition}[theorem]{{\rm P\sc roposition}}

\newtheorem{con}[theorem]{\rm C\sc onjecture}

\newtheorem{para}[theorem]{}

\begin{document}

\title{\vskip-1cm  On surfaces in $\pn{6}$ with no trisecant lines}
\author{ Sandra Di Rocco and Kristian Ranestad}
   \date{}
\maketitle
\thispagestyle{empty}
\begin{center}{\small \sl Dedicated to the memory of F. Serrano}\end{center}
\abstract{Examples of surfaces in $\pn{6}$ with no trisecant lines are constructed. A classification recovering them is given and conjectured to be the complete one.}
\section*{ {\vskip-1cm \Large Introduction}}
The study of varieties embedded in $\pn{N}$ with no trisecant lines
is a very classical problem in algebraic geometry.
The simplest case, i.e. the case of space curves goes back to Castelnuovo.\\
For surfaces the problem has been studied in codimension $2$ and $3$.
In \cite{au} Aure classifies smooth surfaces in $\pn{4}$ with no trisecant through the general point in the space. In \cite{bau} Bauer classifies smooth surfaces in $\pn{5}$ with no  trisecant lines through the general point on the surface.\\
Here we treat the case of smooth surfaces in $\pn{6}$ with no trisecant lines at all. Of course this includes all surfaces which are cut out by quadrics, but there are some more examples.\\
We proceed as follows:
\begin{itemize}
\item In the first two sections we construct examples of surfaces in $\pn{6}$ with no trisecant lines. Standard examples are used to construct new examples, via linkage. In each case we indicate whether the ideal is generated by quadrics or not.
The cases of surfaces containing lines and surfaces with no lines are treated separately.
\item In section \ref{nlist} we give a complete list of surfaces with no proper trisecant lines and no lines on it. 
\item In section \ref{list2} we give a classification of surfaces with lines on the surface but no proper trisecant lines.  These include scrolls, conic bundles, surfaces with an isolated $(-1)$-line and inner projections from $\pn{7}$, i.e.
projections of smooth surfaces in $\pn{7}$ from a smooth point on the surface. 
\end{itemize}
The list of cases produced in sections \ref{nlist} and \ref{list2} cover all the examples constructed in sections \ref{noline} and \ref{wline}. Of course this does not cover all the possibilities.  In fact, surfaces with a finite number of disjoint $(-1)$-lines need not be {\it inner projections} from $\pn{7}$. If in addition every $(-1)$-line meets some other line $L$ of selfintersection $L^2\leq -2$ on the  surface our methods do not apply. However, checking the cases with less than one hundred $(-1)$-lines give no new examples, so we conjecture that the list is in fact complete.\\
The main results of this work are summarized in the following:\\
\newpage
{\rm M{\sc ain  theorem.}} {\it Let $S$ be a smooth surface embedded in $\pn{6}$ with no trisecant lines.
Unless $S$ has a finite number of disjoint $(-1)$-lines, and each one meets some other line $L$ on the surface with $L^2\leq -2$, and $S$ is not an inner projection from $\pn{7}$ the surface belongs to the following list:}
\begin{center}
\begin{tabular}{|c|c|c|c|} \hline
surface &degree&linear system& example \\ \hline\hline
${\bf P}^2$&1,4&${\cal O}_{{\bf P}^2}(1)$, ${\cal O}_{{\bf P}^2}(2)$& \ref{nl1} \\\hline
Rational scrolls&2,3,4,5&linearly normal&  \\\hline
Elliptic scrolls&7&linearly normal& \ref{ell} \\\hline
Rational&4,5,6& anticanonical (Del Pezzo)& \\\hline
$Bl_7({\bf P}^2)$&8&$6l-\sum 2E_i$&\ref{nl2}\\\hline
Rational&6& conic bundle& \ref{con1}\\ \hline
Rational&7& conic bundle& \ref{con2}\\ \hline
Rational&8& conic bundle& \ref{con3}\\ \hline
$Bl_{8}({\bf P}^2)$&8&$4l-\sum_1^8 E_i$& \ref{p48}\\ \hline
$K3$&8&Complete intersection (1,2,2,2)&\ref{nl3}\\\hline
$Bl_{9}({\bf P}^1\times {\bf P}^1 )$&9&$3(1,1)-\sum_1^9 E_i$&\ref{p49}\\ \hline
$Bl_{11}({\bf P}^2)$&10&$6l-\sum_1^5 2E_i-\sum_1^6 E_i$&\ref{p410}\\ \hline
$K3$&10&nontrigonal of genus 6&\ref{nl3},\ref{nl7} \\ \hline
$Bl_1(K3)$&11& nontetragonal of genus 7 $p^{\ast}(\ol{\cH})-E$ & \ref{prk3}\\ \hline 
$Bl_{11}({\bf P}^2)$&12&$9l-\sum_1^5 3E_i-\sum_1^6 2E_i$&\ref{nl4} \\ \hline
$Bl_1(K3)$&12&$p^{\ast}(\ol{\cH})-2E$&\ref{nl5}\\ \hline
Elliptic&12& minimal with $p_g=2$&\ref{elliptic}\\ \hline
Abelian&14&$(1,7)$-polarization&\ref{nl8}\\ \hline
general type&16& complete intersection (2,2,2,2)&\ref{nl9}\\ \hline
\end{tabular}
\end{center}

The authors would like to thank the Mittag-Leffler Institute for its support and its warm environment, which made this collaboration possible and most enjoyable.
%%%%%%%%%%%%%%%%%%%%%%%%%%%%%%%%%%%%%%%%%%%%%%%%%%%%%%%%%%%%%%
\section*{\vskip-1cm \Large Notation}
The groundfield is the field of complex numbers ${\bf C}$.
We use standard notation in algebraic geometry, as in \cite{HAr}.
$S$ is always assumed to be a non singular projective surface.\\
By abuse of notation $\cH_S$ will denote the hyperplane section and the line bundle giving the embedding, with no distinction.\\
When $S$ is the blow up of $\ol{S}$ in $n$ points $S$ will be denoted by $Bl_n(\ol{S})$.
%%%%%%%%%%%%%%%%%%%%%%%%%%%%%%%%%%%%%%%%%%%%%%%%%%%%%%%%%%%%%%%%%%%%%
%%%    EXAMPLES
%%%%%%%%%%%%%%%%%%%%%%%%%%%%%%%%%%%%%%%%%%%%%%%%%%%%%%%%%%%%%%%%%
\section{ \Large Construction of surfaces with no lines}\label{noline}

\subsection{ \large Surfaces defined by quadrics}\label{ex1} 
If $S$ is a surface  whose ideal $I_S$ is generated by quadrics then
clearly it cannot have trisecant lines. The first examples in our
list are then: \\

\begin{para}\label{nl1}The Veronese surface in $\pn{5}$, i.e. $S=\pn{2}$ embedded in
$\pn{5}$ by the linear system $|\cO_{\pn{2}}(2)|$.\end{para}
%%%%%%%%%%%%%%%%%%%%%%%%%%%%%%%%%%%%%%%%%%%%%%%%%%%
\begin{para}\label{nl2} Del Pezzo surfaces of degree $8$ \end{para}
Let $S$ be the blow up of $\pn{2}$ in $7$ points embedded in $\pn{6}$ by the
linear system $|-2K_S|$.
The line bundle $\cH_S=-2K_S$ is $2$-very ample and thus embeds $S$ without
trisecant lines. See \cite{sa} for the definition and the proof of the $2$-very ampleness.
One can construct this surface in $\pn{6}$ as the intersection of the cone over a Veronese surfaces in $\pn{5}$ with a quadric hypersurface. Thus the surface is cut out by quadrics, and it has no lines as soon as the quadric does not contain the vertex of the cone.  These surfaces are cut out by 7 quadrics in $\pn{6}$. 

%%%%%%%%%%%%%%%%%%%
\begin{para}\label{nl3} General minimal nontrigonal {\it $K3$-surface of degree $8$ or $10$}.\end{para} Consider a nontrigonal $K3$ surface of degree $8$ in $\pn{5}$.  It is the complete intersection of $3$ quadrics, and the general one has Picard group generated by the hyperplane section so it has no lines.\par
Similarly a nontrigonal $K3$ surface of degree $10$ in $\pn{6}$ is a linear section of the Pl\"ucker embedding of the Grassmannian $Gr(2,5)$ intersected with a quadric hypersurface.  Again the surface is cut out by quadrics, in fact 6 quadrics,  and the general one has Picard group generated by the hyperplane section, so it has no lines. 
%%%%%%%%%%%%%%%%%%%%%%%%%%%%%%%%%%%%%%%%%%%%%%%%%%%%%%%%%
\begin{para}\label{nl7}Minimal tetragonal $K3$-surfaces of degree $10$\end{para}
Let $V$ be the cubic cone over $\pn{1}\times\pn{2}$, and consider two general quadrics $\Qs{1},\Qs{2}$, containing a quadric surface $S_2$ in $\pn{3}$. The complete intersection $V\cap \Qs{1}\cap\Qs{2}=S\cup S_2$ produces a smooth surface $S$ of degree $10$ in $\pn{6}$. A trisecant line of $S$ would be a line in $S_2$.  But the intersection curve $S\cap S_2$ is an elliptic quartic curve, of type $(2,2)$ on $S_2$, so there is no trisecant. Again using adjunction one can see that $p_g(S)=1$ and $K_S\cdot \cH_S=0$ and that the ruling $\cF$ of $V$ restricted to $S$ form a quartic elliptic curve on $S$. It follows that $S$ is a $K3$ surface of degree $10$ in $\pn{6}$, with a pencil of elliptic quartic curves on it.  Furthermore it is cut out by 6 quadrics.\\
One can also construct $S$ by linkage from the rational surface $Bl_7(\pn{2})$
 embedded in $\pn{5}$ by the line bundle $4l-2E-\sum_1^6E_i$. 
  These surfaces are cut out by 4 quadrics in $\pn{5}$ ( see \ref{con1}), and in a complete intersection $(2,2,2,2)$ in $\pn{6}$ they are linked to $K3$-surfaces of the above type. 
 
%%%%%%%%%%%%%%%%%%%%%%%%%%%%%%%%%%%%%%%
%%%%%%%%%%%%%%%%%%%%%%%%%%%%%%%%%%%%%%%%%%%%%%%%%%%%%%%%%%%%%%%%%%%%%%%%%%  
\begin{para}\label{elliptic} Two families of elliptic surfaces of degree $12$ \end{para}
%%%%%%%%%%%%%%%%%%%%%%%%%%%%%%%%%%%%%%%%%%%%%%%%%%%%%%%%%%%%%%%%%%%%%
%\begin{para}\label{ellipticII}Elliptic surface of degree $12$ (II)\end{para}
Let $V$ be a rational normal $4$-fold scroll of degree $3$ in $\pn{6}$ and consider $S=V\cap{\cal Q}_1\cap{\cal Q}_2$, where 
${\cal Q}_1,{\cal Q}_2$  are two general quadrics, which do not have any common point in the singular locus of $V$.  Note that this is possible only if the cubic 4-fold has vertex a point or a line.  This gives us two separate cases.  In both cases the intersection of $V$ with the two quadrics is a smooth surface. Now $K_{\ol V}=-4\cH+\cF$, where $\ol V$ is the $\pn{3}$-bundle over $\pn{1}$ which is mapped to $V$ by $\cH$, and $\cF$ is a member of the ruling.
Then $K_S=\cF|_S$ gives a fibration of elliptic quartic curves without multiple fibers onto $\pn{1}$.\\
$S$ does not have proper trisecant lines since it is cut out by quadrics.
Moreover $S$ cannot contain lines. In fact if there were a line $L\subset S$, then $L\subset V$ and $L\subset {\cal Q}_i$. The biggest component of the Fano variety of lines in $V$ is the 5-dimensional component of lines in the pencil of $\pn{3}$s.  Since a line imposes three conditions on a quadric, 
 we see that in order for $L$ to be on two quadrics we need at least a $6$-dimensional family of lines on $V$.  Thus we have constructed two families of elliptic surfaces of degree 12, one on the cubic 4-fold cone with vertex a point, and one on the cubic 4-fold cone with vertex a line.  In the first case any two canonical curves span $\pn{6}$, while in the other case any two canonical curves span a $\pn{5}$.   These surfaces are cut out by 5 quadrics in $\pn{6}$. 
 %%%%%% 
\begin{para}\label{nl9} Complete intersections of $4$ quadric hypersurfaces in $\pn{6}$.\end{para}
This is a degree $16$ surface of general type. Since the Fano variety of lines in a quadric has codimension $3$ in the Grassmannian of lines in $\pn{6}$ and this Grassmannian has dimension $10$, there
are no lines in a general complete intersection of $4$ quadrics.  A similar argument is made more precise in \ref{nl5} below.
%%%%%%%%%%%%%%%%%%%%%%%%%%%%%%%%%%%%%%%%%%%%%%%%%%%%%%%%%%%%%%%%%%%%%%%%%%%%%%%
\subsection{\large Non minimal $K3$-surfaces of degree 12}\label{nl5}
  Consider $4$ general quadrics $\cQ_1,...,\cQ_4\sub\pn{6}$ containing the Veronese surface  $S_4\sub\pn{5}$ and let $S$ be the residual surface of degree 12 in $\pn{6}$.  Let $V$ be the complete intersection of a general subset of three of the quadrics $\cQ_1,...,\cQ_4$.  
By adjunction we see that $C=S\cap S_4=\cH_V|_{S_4}-K_{S_4}= 5l$, where $l$ is the generator of $Pic(S_4)$. Then $K_S=(\cH_V-S_4)_S$ and $K_S\cdot \cH_S=((\cH_V+S_4)_S)\cdot \cH_S=12-10=2$.
The exact sequence 
\begin{equation}\label{E1}
0\to\cO_V(\cH_V-S_4)\to\cO_V(\cH)\to\cO_{S_4}(\cH)\to 0
\end{equation} and the fact that $h^0(\cO_{S_4}(\cH))=6$ gives  $h^0(\cO_V(\cH_V-S_4))=1$ and $h^1(\cO_V((\cH_V-S_4))=h^2(\cO_V(((\cH_V-S_4))=0$. Plugging those values in the long exact sequence of
\begin{equation}\label{E2}
0\to\cO_V(-\cH)\to\cO_V(\cH_V-S_4)\to\cO_S(K_S)\to 0
\end{equation}
we get $p_g(S)=1$, $q=0$ and thus $\chi(S)=2$.  Thus the canonical curve is a $(-1)$- conic section or two disjoint $(-1)$-lines.  It is the residual to the intersection $C$ of $S$ with $S_4$ in a hyperplane section of $S$, so each component intersect the curve $C$ in at least two points.  If the canonical curve is two disjoint lines these would be secants to $S_4$. But a secant line to $S_4$ is a secant line to a unique conic section on $S_4$. Any quadric which contains a secant line must therefore contain the plane of this conic section.  Therefore $S$ contains no secant line to $S_4$ as soon as $S$ is irreducible. We conclude that the canonical curve is a $(-1)$-conic section.  The surface $S$ is the blow up of a minimal $K3$-surface and $K_S^2=-1$.
We now show that $S$ contains no lines.  First assume that the line does not intersect $S_4$.
Let $G$ be the Grassmannian of $4$-dimensional subspaces of quadrics in $\pn{6}$ which contain $S_4$.  Consider the incidence variety $$I=\{(L,Q_4)\in Gr(2,7)\times G|L\sub\cap_{\cQ\in Q_4} \cQ, L\cap S_4=\emptyset\}$$
and let $p:I\to Gr(2,7)$ and $q: I\to G$ the two projections.
Since the lines in $\pn{6}$ form a $10$-dimensional family and $L\sub \Qs{i}$ is a codimension $3$ condition it is clear that codim$(q(p^{-1}(L)))=12$ and that codim$(q(p^{-1}(Gr(2,7)))\geq 12-10=2$. This means that we can choose $Q_4$ general enough so that there is no line $L\sub\cap_{\cQ\in Q_4} \cQ$ which is disjoint from $S_4$.\\
Similarly consider $I=\{(L,Q_4)\in Gr(2,7)\times G|L\sub\cap_{\cQ\in Q_4} \cQ,\, L\cap S_4\not=\emptyset\}$ and $p,q$ as before where now $p:I\to \ol{G}=
\{L\in Gr(2,7)|L\cap S_4\not=\emptyset\}$. The intersection of a line in $\pn{6}$ with a surface imposes $3$ conditions
 so that dim$\ol{G}=7$.  Since a line through a point on $S_4$ imposes 2 conditions on quadrics through $S_4$ we have codim$(q(p^{-1}(l)))=8$ in this case, it follows that codim$(q(p^{-1}(\ol{G})))\geq 8-7=1$ so we can assume that there are no lines on $S$ intersecting $S_4$ in one point.
We are left to examine the case when $L$ meets $S_4$ in at least 2 points, i.e. when $L$ is a secant line to the Veronese surface $S_4$.  But as above this is impossible as long as $S$ is irreducible.\\
 Thus $S$ has no lines on it and it is the blow up of a $K3$ surface in one point and $\cH_S=p^*(\ol{\cH}-2E)$, where $p:S\to\ol{S}$ is the blow up map and $\ol{\cH}$ is a line bundle on $\ol{S}$ of degree $16$.  Again $S$ has no trisecant line because such a line would necessarily be a line in the Veronese surface.  Notice that the conic sections on the Veronese surface each meet $S$ in 5 points.  In fact the Veronese surface is the union of the 5-secant conic sections to $S$ and is therefore contained in any quadric which contains $S$.  It is straightforward to check that the surfaces $S$ are cut out by $4$ quadrics and $3$ cubics in $\pn{6}$.

%%%%%%%%%%%%%%%%%%%%%%%%%%%%%%%%%%%%%%%%%%%%%%%%%%%%%%%%%%%%%%%%%%%%%%%%%%%%%
\subsection{\large A family of rational surfaces of degree 12}\label{nl4}
Let $S=Bl_{11}(\pn{2})$ be polarized by the line bundle $H=9l-\sum_1^5 3E_i-\sum_6^{11}2E_j$. Assume that the $11$ points blown up are in general position. More precisely we require that the following linear systems are empty for all possible sets of distinct indices:
\begin{itemize}
\item $|E_i-E_j|$
\item $|l-E_i-E_j-E_k|$
\item $|2l-\sum_{k=1}^6 E_{i_k}|$
\item $|3l-\sum_{k=1}^{10}E_{i_k}|$ and $|3l-2E_{i_0}-\sum_{k=1}^7E_{i_k}|$
\item $|4l-2\sum_{k=1}^2E_{i_k}-\sum_{k=3}^{11}E_{i_k}|$ and $|4l-2\sum_{k=1}^3E_{i_k}-\sum_{k=4}^{10}E_{i_k}|$
\item $|5l-2\sum_{k=1}^5E_{i_k}-\sum_{k=6}^{11}E_{i_k}|$ 
\item $|6l-2\sum_{k=1}^{8}E_{i_k}-\sum_{k=7}^{11}E_{i_k}|$, $|6l-3E_{i_0}-2\sum_{k=1}^{5}E_{i_k}-\sum_{k=6}^{10}E_{i_k}|$,\\ $|6l-2\sum_{k=1}^{9}E_{i_k}-E_{i_{10}}|$
\end{itemize}
\begin{lemma} $H$ is a very ample line bundle on $S$.
\end{lemma}
\proof $H$ is shown to be very ample in \cite{li}. We report here a different short proof which relies on 
\begin{lemma} {\rm (Alexander)\cite[Lemma 0.15]{kr}.}  If $H$ has a decomposition $$H=C+D,$$
where $C$ and $D$ are curves on $S$, such that ${\rm dim}|C|\geq 1$, and if the restriction maps $H^0(\cO_S(H))\to H^0(\cO_D(H))$ and $H^0(\cO_S(H))\to H^0(\cO_C(H))$ are surjective, and $|H|$ restricts to a very ample linear system on $D$ and on every $C$ in $|C|$, then $|H|$ is very ample on $S$. \end{lemma} \
Consider the reducible hyperplane section:
$$H=D_1+D_2=(3l-\sum_1^8E_i)+(6l-2\sum_1^5E_i-\sum_6^8E_j-2\sum_9^{11}E_k)$$
Then $D_1$ is embedded as a degree 6 elliptic curve and $D_2$ as a sextic curve of genus $2$. Moreover a general element of $|D_2|$ is irreducible by the choice of points in general position and all the elements in the pencil $|D_1|$ are irreducible. It follows that $H_{D_1}$ and $H_{D_2}$ are very ample. The fact that both the maps $H^0(S,H)\to H^0(D_1,H_{D_1})$ and $H^0(S,H)\to H^0(D_2,H_{D_2})$ are surjective concludes the argument.
\qed 
\begin{lemma} There are no lines on $S$.\end{lemma}
\proof Assume $L=al-\sum_1^{5}a_i E_i-\sum_6^{11}b_j E_j$ is a line on $S$.
 Then looking at the intersection of $L$ with the twisted cubics $E_i$, $i=1,...,5$, and the intersection of $L$ with the conics $E_j$, $j=6,...,11$ we derive the bounds $0\leq a_i\leq 2$ and $0\leq b_i\leq 2$.
This implies that $1=H\cdot L=9a-\sum_1^5 3a_i-\sum_6^{11}2b_i\geq 9a-30-24$, i.e. $a\leq 4$. The only numerical possibilities are:
\begin{itemize}
\item $L=E_i-E_j$;
\item $L=3l-2E_{i_1}-\sum_{k=2}^5E_{i_k}-\sum_{k=6}^9E_{i_k}\quad \{i_1,\dots ,i_5\}=\{1,\dots ,5\}$,\\ $6\leq i_6<\dots <i_9\leq 11$;
\item $L=2l-\sum_1^5E_i-E_j$, $6\leq j\leq 11$;
\item $L=l-E_i-E_j-E_k\quad 1\leq i<j\leq 5<k\leq 11$.
\end{itemize}
But those are empty linear systems by the general position hypothesis.
\qed
\begin{proposition}$S$ is a rational surface in $\pn{6}$ with no trisecant lines.
\end{proposition}
\proof Consider the reducible hyperplane section $H=\Fd+\Fu$ where
\begin{itemize}
\item $\Fd=-K_S+E_i+E_j$ $i,j=6\dots,11$, i.e. an elliptic quartic curve.
\item $\Fu=H-\Fd=6l-2\sum_1^5E_k-\sum_6^{11}E_k-E_i-E_j$, i.e. a curve of degree $8$ of genus $3$.
\end{itemize}
Fix a curve $\Fd$. Any trisecant line $L$, would together with $\Fd$ span a hyperplane, so there is some reducible hyperplane section $H=\Fd+\Fu$ for which $L$ is a trisecant. Since neither $\Fu$ nor $\Fd$ have trisecants, $L$ must in fact intersect both of these curves, so $\Fd$ and $L$ spans at most a $\pn{4}$.
This means that we can always find a curve $C\in |\Fu|$ passing through the points in $L\cap\Fd$, and such that $L$ is a trisecant to $C\cup \Fd$.
If $C$ is irreducible this is impossible since $C$ has no trisecant lines\\
Assume it is reducible and write $C=A+B$, were $A$ and $B$ are irreducible with $deg(A)\leq deg(B)$. Then  the following cases could occur:
\begin{itemize}
\item[(a)] $A$ is a plane conic;
\item[(b)] $A$ is a plane cubic or a twisted cubic;
\item[(c)] $A$ is a quartic curve;
\end{itemize}
Let $A=\al l-\sum_1^5\as{i}E_i-\sum_6^{11}\be{j}E_j$.\\
Assume $\al=0$. If $A=E_k$ for $k\in\{1,...,5\}$ then $B=\Fd -E_k$, which is impossible by the general position hypothesis. If $A=E_k$ for $k\in\{6,...,11\}\setminus\{i,j\}$ then $B=\Fu -E_k$, this possibility will be analyzed more closely below.\\
Assume now $\al>0$, i.e. $A\neq E_i$, then by intersection properties and the assumption that $A$ and $B$ are effective divisors, $0\leq\as{i}\leq 2$, $0\leq\be{j}\leq 1$ and $1\leq\al\leq 6$. \\
\noindent (a) Going over the possibilities for $\al,\as{i}$ and $\be{j}$ gives no result by the general position hypothesis.\\
\noindent (b) Examining the possible choices for $\al,\as{i}$ and $\be{j}$
 we get: 
\begin{itemize}
\item $A=l-E_m-E_n, 1\leq m<n\leq 5$ residual to $B=5l-2\sum_{k=1}^5E_k+E_m+E_n-\sum_{k=6}^{11} E_k-E_i-E_j$
\item $A=2l-\sum_1^5E_k$ residual to $B=4l-\sum_1^{11}E_k-E_i-E_j$
\item $3l-\sum_{k=1}^{11}E_k+E_m, 1\leq m\leq 5$ and $3l-\sum_1^5E_k-E_m-E_i-E_j$
\item $4l-2\sum_{k=1}^{11}E_k-E_m-E_n, 1\leq m<n\leq 5$  
\item $5l-2\sum_1^5E_i-\sum_6^{11}E_j$ 
\end{itemize}
In the first two cases the residual curve $B$ does not exist by the general position hypothesis.  Likewise the curve $A$ does not exist in the last cases.\\
\noindent(c) Similar computations leads to
\begin{itemize}
\item $A=l-E_m-E_n$, $1\leq m\leq 5$ and $6\leq n\leq 11$, whose residual curve does not exist.
\item $A=2l-\sum_1^5E_k+E_m-E_n$, $1\leq m\leq 5<n\leq 11$, whose residual curve does not exist.
\item $A=F_{k,l}$, $B=F_{m,n}$ with $\{i,j,k,l,m,n\}=\{6,\dots,11\}$
\end{itemize}
We are left with the cases $C=E_k+(6l-2\sum_1^{5}E_m-\sum_6^{11}E_n-E_i-E_j-E_k)$ or $C=F_{k,l}+F_{m,n}$.\\
Notice that in both cases the projective spaces spanned by the two components, $<A>,<B>$, intersect in a line $L=<A>\cap<B>$. Moreover neither $A$ nor $B$ admits trisecant lines and $A\cap B=2$.
It follows that $A\cap <B>=A\cap B=B\cap <A>$. Any trisecant line $L$ to $C$, must meet $A$ ( or respectively $B$) in two points and thus it is contained in $<A>$, which implies $L\cap B\sub A\cap B$. But this means that $L$ is a trisecant line for $A$, which is impossible. \qed
 These surfaces are cut out by 3 quadrics and 4 cubics in $\pn{6}$. 

\subsection{ \large Abelian surfaces}\label{nl8}
Recently Fukuma \cite{ab1}, Bauer and Szemberg \cite{ab2} have proved that the general $(1,7)$-polarized abelian surface in $\pn{6}$ does not have any trisecant lines.  The argument uses a generalization of Reider's criterium to higher order embeddings. These surfaces are not contained in any quadrics. 

%%%%%%%%%%%%%%%%%%%%%%%%%%%%%%%%%%%%%%%%%%%%%%%%%%%%%%%%%%%%%%%%%%%%%%%%%%%%%
%%%%%%%%%%%%%%    SURFACES WITH LINES %%%%%%%%%%%%%%%%%%%%%%%%%%%%%%%%%%%%%%%%%
\section{ \Large Construction of surfaces with lines}\label{wline}
In this section we will construct examples of surfaces with no proper trisecant lines but with lines on it. In all the cases below the surfaces are cut out by quadrics, so naturally there are no proper trisecant lines.  In fact, we do not know of
any surface with lines on it, but with no proper trisecant which is not cut out by quadrics.  

Natural examples of surfaces  with at least a one dimensional family of lines are 
given by surfaces of minimal degree, i.e. the rational normal scrolls of degree $N-1$ in $\pn{N},3\leq N\leq 6$. 
Similarly the Del Pezzo surfaces of degree 4, 5 and 6 form natural families of surfaces cut out by quadrics.  They have a finite number of lines on them.  We construct other rational and nonrational surfaces.
 %%%%%%%%%%%%%%%%%%%%%%%%%%%%%%%%%%%%%%%%%%%%%%%%%%%%%%%%%%%%%%%%%%%%%
\subsection{\large Rational surfaces}
%%%%%%%%%%%%%%%%%%%%%%%%%%%%%%%%%%%%%%%%%%%%%%%%%%%%%%%%
\begin{para}\label{con1}Conic bundles of degree $6$\end{para} Consider a cubic 3-fold scroll $V\sub\pn{5}$ and let $S=V\cap \cQ$ be the surface of degree $6$ given by the intersection with a general quadric hypersurface.  $S$ is smooth as soon as $\cQ$ avoids the vertex of $V$.  Therefore we have two cases, when $V$ is smooth and when $V$ is the cone over a smooth cubic surface scroll in $\pn{4}$.\\
Let $\ol{V}$ be the $\pn{2}$-scroll which is mapped to $V$ by $\cH$, and let $\cF$ be a member of the ruling. By abuse of notation we denote the pullback of $S$ to $\ol{V}$ by $S$, it is isomorphic anyway.  Then by adjunction 
$$K_S=(-3\cH+\cF+2\cH)|_S=-\cH|_S+F|_S$$
and thus $(K_S)^2=2$ and $K_S\cdot\cH_S=-4$.  Furthermore $p_g(S)=q(S)=0$, so $S$ is rational.
The ruling of the scroll define a conic bundle structure on $S$, and there are 6 singular fibers, i.e. 12 $(-1)$-lines in the fibers since  $(K_S)^2=2$. These surfaces are cut out by 4 quadrics in $\pn{5}$.  %%%%%%%%%%%%%%%%%%%%%%%%%%%%%%%%%%%%%%%%%%%%%%%%%%%%%%%%%%%%%%%%% 
\begin{para}\label{con2}Conic bundles of degree $7$\end{para} 
Consider a rational normal 3-fold scroll  $V\sub\pn{6}$ of degree 4, and let $\cQ$ be a general quadric hypersurface containing a member of the ruling $\cF$. Then the complete intersection $V\cap\cQ=\cF\cup S$, where $S$ is a surface of  degree $7$ in $\pn{6}$.  As soon as $V$ is smooth, $S$ is also smooth.  In fact in this case $S$ would be singular if $V$ was singular.\
By adjunction 
$$K_S=(-3\cH+\cF+2\cH)|_S=-\cH|_S+F|_S,$$
where $\cH$ is a hyperplane section.
Thus $(K_S)^2=3$ and $K_S\cdot\cH_S=-5$. 
Like in the previous case the ruling of $V$ define a conic bundle structure on $S$ with 5 singular fibers, i.e. 10 $(-1)$-lines altogether in the fibers. If $S$ had a trisecant line, $L$, then $L$ would be a line in $\cF$ intersecting the curve $C=S\cap\cF$ in three points. In this case $C$ is conic section so there is no trisecant.  In fact one can also show that $S$ is cut out by 8 quadrics.
%%%%%%%%%%%%%%%%%%%%%%%%%%%%%%%%%%%%%%%%%%%%%%%%%%%%%%%%%%%%%%%%%%%%%%%%%%
 \begin{para}\label{con3}Conic bundles of degree $8$\end{para} 
Let $V$ be a rational normal $3$-fold scroll of degree $4$ in $\pn{6}$, and let $\cQ$ be a general quadric hypersurface not containing any singular point on $V$.  Then $V$ is smooth or is a cone with vertex a point, and $S=V\cap\cQ$ is a smooth surface of degree $8$ in $\pn{6}$.  We get two cases like in \ref{con1}.  Proceeding with notation like in that  case  we get $K_S=2\cF_S-\cH_S$ and $K_S^2=0$.  Thus we get conic bundles with sectional genus $3$ with $8$ singular fibers, i.e. 16 $(-1)$-lines in fibers.  These surfaces are cut out by $7$ quadrics. 
 
%%%%%%%%%%%%%%%%%%%%%%%%%%%%%%%%%%%%%%%%%%%%%%%%%%%%%%%%%%%%%%%%%%%%%%%%%%%%%
\begin{para}\label{p48}A family of surfaces of degree $8$\end{para} Let $V$ the a cone over $\pn{1}\times\pn{2}$ as in the previous example  and let 
$\Qs{1},\Qs{2}$ be two quadrics containing a $F=\pn{3}$ of the ruling, in particular they pass through the vertex on the cone. Then $V\cap\Qs{1}\cap\Qs{2}=S\cup F$, where $S$ is a smooth surface of degree $8$ in $\pn{6}$. In fact let $E$ be the exceptional divisor in the blow up $\ol{V}$ of $V$ at the vertex, and let $\ol{S}$ be the strict transform of $S$. With notation as above $\ol{S}=(2\cH-\cF-E)\cap (2\cH-\cF-E)$ on $\ol{V}$.  The exceptional divisor $E$ is isomorphic to $\pn{1}\times\pn{2}$ and 
$E|_E=-\cH$ via this isomorphism. Furthermore the self intersection of $E|_{\ol{S}}$ is 
\[\begin{array}{ll}[(2\cH_V-\cF-E)^2]\cdot E\cdot E&=[4\cH^2-4\cH F+2\cF E-a\cH E+E^2]\cdot E\\
&=2H^2\cdot\cF-E^3=2-3=-1
\end{array}  \]

By adjunction $K_S=-\cF_S$ and thus $|-K_S|$ is a pencil of elliptic curves with one base point at the vertex of $V$. It follows that $S$ is rational of degree $8$, sectional genus $3$ and $K_S^2=1$. The adjunction $|\cH+K_S|$ maps the surface birationally to $\pn{2}$, so the the surface is $Bl_{8}({\bf P}^2)$ with $\cH=4l-\sum_1^8 E_i$. These surfaces are cut out by 6 quadrics in $\pn{6}$.

%%%%%%%%%%%%%%%%%%%%%%%%%%%%%%%%%%%%%%%%%%%%%%%%%%%%%%%%%%%
\begin{para}\label{p49}Two families of surfaces of degree $9$\end{para}
Let $\ol{V}$  be a $\pn{3}$-bundle of degree $3$ over $\pn{1}$ with ruling  $\cF$ and let $V$ be its image rational normal $4$-fold of degree $3$ in $\pn{6}$ under the map defined by $\cH)$ as in \ref{elliptic}. Let $\Qs{1},\Qs{2}$ be general quadrics with no common point in the vertex of $V$ which contain a smooth cubic surface $S_3=V\cap\pn{4}$, for some general $\pn{4}\sub\pn{6}$. Then $V\cap\Qs{1}\cap\Qs{2}=S_3\cup S$, and $S$ is smooth as soon as the vertex of $V$ is a line or a point.  The curve of intersection $C=S_3\cap S$ is then a curve of degree $6$ represented in $S_3$ by the divisor $2\cH_{S_3}$. Thus $S$ cannot have trisecants since $C$ has no trisecant lines.\\
Since $S_3$ meets each ruling of $V$ in a line, $S$ which is linked to $S_3$ in 2 quadrics meets each ruling of $V$ in a twisted cubic curve.  Therefore $S$ is rational. Furthermore,
$K_{\ol{V}}=-4\cH+\cF$, so by adjunction $K_S=(-4\cH+\cF+4\cH)|_S-S_3|_S=-C+\cF_S$.  Thus 
 $K_S\cdot\cH_S=(\cF_S-C)\cdot \cH_S=-3$, and $S$ has sectional genus $4$. Notice that, by adjunction, $C\cdot \cF_S=2$. Let $D=\cH-C$, then $D$ has degree 3 and moves in a pencil, so it must be a twisted cubic curve, with $D^2=0$.  The genus formula implies that $C\cdot D=3$, and so $C^2=3$.  Therefore $K_S^2=(\cF_S-C)^2=-1$.  
>From the two types of cubic scrolls $V$ we get two types of surfaces $S$. Both are rational surfaces of degree $9$, sectional genus $4$ and $K_S^2=-1$.  In both case $S$ has two pencils of twisted cubic curves, but in one case any two of the curves in a pencil span $\pn{6}$, while in the other case any two of them span a $\pn{5}$.  The adjunction $|\cH+K_S|$ maps the surface birationally to a smooth quadric surface in $\pn{3}$, so the surface is $Bl_{9}({\bf P}^1\times {\bf P}^1)$ with $\cH=3(1,1)-\sum_1^9 E_i$.  The two families correspond to the cases when the 9 points on the quadric is the complete intersection of two rational quartic curves (of type $(3,1)$ and $(1,3)$ respectively) or not. These surfaces are cut out by 6 quadrics in $\pn{6}$.

%%%%%%%%%%%%%%%%%%%%%%%%%%%%%%%%%%%%%%%%%%%%%%%%%%%%%%%%%%%%%%%%%%%%%%%%%%%%%
\begin{para}\label{p410}A family of surfaces of degree $10$\end{para}
Consider the Del Pezzo surface $S_6$ of degree $6$ in $\pn{6}$, and $4$ general quadrics containing it,  $\Qs{1},\Qs{2},\Qs{3},\Qs{4}$. Then $\Qs{1}\cap\Qs{2}\cap\Qs{3}\cap\Qs{4}=S\cup S_6$, where $S$ is a smooth surface of degree $10$ in $\pn{6}$. The exact sequences (\ref{E1}), (\ref{E2}) of \ref{nl5} applied in this case show that $S$ is rational.  The adjunction formula gives $K_S\cdot \cH_S=(\cH-S_6)_S\cdot\cH_S=-2$, so the sectional genus is $5$. The intersection curve $C=S\cap{S_6}=2\cH_{S_6}$ on $S_6$,  which means that $S$ has no trisecant lines.  These surfaces are cut out by 5 quadrics in $\pn{6}$. The adjoints of the surfaces in \ref{nl4} are of this type.
%%%%%%%%%%%%%%%%%%%%%%%%%%%%%%%%%%%%%%%%%%%%%%%%%%%%%%%%%%%%%%%%%%%%%%%%%%%%%%
\subsection{\large Non rational surfaces}
%%%%%%%%%%%%%%%%%%%%%%%%%%%%%%%%%%%%%%%%%%%%%%%%%%%%%%%%%%%%%%%%%%%%%%%%%%%%%
There are also non rational surfaces without proper trisecant lines.  The first examples are
\begin{para}\label{ell} Elliptic scrolls\end{para}
The elliptic normal scrolls of degree 7 for which the minimal self intersection of a section is $1$, are cut out by $7$ quadrics (cf. \cite{hks}).\\
%%%%%%%%%%%%%%%%%%%%%%%%%%%%%%%%%%%%%%%%%%%%%%%%%%%%%%%%%%%%%%%%%%%%%%%%%%%%%%
Finally there are
%%%%%%%%%%%%%%%%%%%%%%%%%%%%%%%%%%%%%%%%%%%%%%%%%%%%%%%%%%%%%%%%%%%%%%%%%%%%%
\begin{para}\label{prk3}A family of non minimal $K3$-surfaces.\end{para}
Consider an inner projection of a general nontrigonal and nontetragonal $K3$-surface $\ol{S}$ of degree $12$ in $P^7$ cf. \cite{Mu}. The surface $S$ is the projection from a point $p\in \ol{S}$, $\pi_p:\ol{S} \to S$. Then $S$ is a $K3$-surface of degree $11$ in $\pn{6}$ with one line, i.e. the exceptional line over $p$.  Any trisecant of $S$ will come from a trisecant of $\ol{S}$ or from a four secant $\pn{2}$ to $\ol{S}$ through $p$.  But a normally embedded $K3$-surface with a trisecant is trigonal and with a $4$-secant plane is tetragonal which is avoided by assumption. So $S$ has no trisecant. These surfaces are cut out by $5$ quadrics in $\pn{6}$. 

%%%%%%%%%%%%%%%%%%%%%%%%%%%%%%%%%%%%%%%%%%%%%%%%%%%%%%%%%%%%%%%%%%%%%
%%%%%%%%%%%%%%%%%%%%%%%%%%%%%%%%%%%%%%%%%%%%%%%%%%%%%%%%%%%%%
%LIST
%%%%%%%%%%%%%%%%%%%%%%%%%%%%%%%%%%%%%%%%%%%%%%%
\section{ \Large   Complete list of surfaces with no lines}\label{nlist}
Throughout this section we will assume that $S$ is a surface embedded in $\pn{6}$ by the line bundle $\cH$, with no lines on it and no trisecant lines.  $C$ will denote the general smooth hyperplane section of $S$.\\
With this hypothesis the invariants of the surfaces will give zero in the two formulas of Le Barz in \cite{leb}. Let
\begin{itemize}
\item $n=$degree$(S)$; 
\item      $k=K_S^2$;
\item      $c=c_2(S)$; 
\item      $e=K_S\cdot\cH$;
\end{itemize}
Then the formula of Le Barz for the number of trisecant lines meeting a fixed $\pn{4}\sub\pn{6}$ is:
\begin{equation}\label{A}
 D_3=2 n^3  - 42 n^2  + 196 n-k(3n-28)+c(3n-20) - e(18n - 132)\quad (=0)
\end{equation}
and the formula for the number of lines in $\pn{6}$ which are tangential trisecants, i.e. tangent lines that meet the surface in a scheme of length at least 3, is:
\begin{equation}\label{B}
 T_3=6 n^2  - 84 n+k(n-28) -c(n-20)+ e(4n-84)\quad (=0)
\end{equation}
We set these equal to 0.
Next we use the Castelnuovo bound for an irreducible curve of degree $n$ in $\pn{N}$ \cite{ACGH}:
$$p(N)=[\frac{n-2}{N-1}](n-N-([\frac{n-2}{N-1}]-1)\frac{N-1}{2})$$
where $[x]$ means the greatest integer $\leq x$, and refined versions of it given by the following two theorems of Harris and Ciliberto:
\begin{theorem}\label{HA}\cite{ha}, \cite[3.4]{ci}
Let $p_1=\frac{n^2}{10}-\frac{n}{2}$ and let $C$ be a reduced, irreducible curve in $\pn{5}$ of degree $n$ and genus $g$. Then
\begin{itemize}
\item[(a)] If $g(C)> p_1$  then $C$ lies on a surface of minimal degree in $\pn{5}$;
\item[(b)] if $g=p_1$ and $n\geq 13$ then $C$ lies on a surface of degree $\leq r$.
\end{itemize}
\end{theorem}
\begin{theorem}\cite[Th. 3.7]{ci}\label{p3}
Let $p_3(n,r)=\frac{n^2}{2(r-1)+3}+O(n)$ and $p_2(n,r)=\frac{n^2}{2(r-1)+4}+O(n)$, where $0\leq O(n)\leq 1$.\\
Let $C$ be a reduced, irreducible, non degenerate curve in  $\pn{6}$, $r\geq 6$, of degree $n$ and genus $p$. Then
\begin{itemize}
\item  if $n>(r+1)$ and $p>p_3(n,r)$ then $C$ lies on a surface of degree $\leq r$;
\item if $2r+3\leq n\leq 5r+2$ and $p>p_2(n,r)$ then $C$ lies on a surface of degree $\leq r$.
\end{itemize}  
\end{theorem}

\subsection{ \large The cases with $n\leq 11$}
Notice that surfaces in $\pn{3}$ and in $\pn{4}$ have trisecant lines or contain lines.  The only surfaces in $\pn{5}$ which do not contain lines or have  trisecants are the Veronese surfaces and the general complete intersections $(2,2,2)$, i.e. the general nontrigonal $K3$-surfaces of degree 8.\\
 Then assuming $N\geq 5$ 
and imposing the following conditions :
\begin{itemize}
\item (\ref{A}) and (\ref{B});
\item $g\leq p(5)$ ($S$ spans $\pn{6}$);
\item $c+k=12\cdot$integer;
\item $k\leq 3c$ (Miyaoka) and $kn\leq e^2$ ( Hodge index Theorem)
\end{itemize}
numerical computations give a list of possible invariants for $n\leq 11$:
\begin{center}
\begin{tabular}{|c||c|c|c|c|} \hline
&$n $ & $e$& $k$& $c$\\ \hline\hline
(1)&$4$ &$ -6$ &$9$& $3$\\ \hline
(2)&$ 8$& $ -4$& $2$& $10$\\ \hline
(3)&$ 8$&$0$&$ 0$&$ 24$ \\\hline
(4)&$10$ &$0$&$ 0$&$ 24$ \\\hline
\end{tabular}
\end{center}
The examples \ref{nl1}, \ref{nl2}, \ref{nl3}, \ref{nl7} in section 1 have these invariants, and it is easy to see that these are the only ones.  In the first two cases any smooth surface in the family would have no lines, but for $K3$-surfaces it is easy to construct degenerations to smooth surfaces with one or several lines.  These lines would be $(-2)$-lines on the $K3$-surface.
Therefore the four cases are:

\begin{itemize}
\item $(S,\cH)=(\pn{2},\Os{\pn{2}}(2)$ is a Veronese surface in $\pn{5}$;
\item $(S,\cH)=(Bl_7(\pn{2}),-2K_S)$ is a Del Pezzo surface in $\pn{6}$;
\item $S$ is a general nontrigonal $K3$ surface of degree $8$ in $\pn{5}$;
\item $S$ is a general tetragonal or nontetragonal $K3$ surface of degree $10$ in $\pn{6}$. 
\end{itemize} 
\subsection{ \large When $n\geq 12$}
Eliminating $c$ from (\ref{A}) and (\ref{B}) gives 
$$k=\frac{[n^4-32n^3+332n^2-1120n]-e(3n^2-80+480)}{8n}$$
By \ref{HA} if $n\geq 11$ then we get for the general hyperplane section $C$ of $S$ that $g(C)\leq p_1$ or $C$ lies on a surface of minimal degree in $\pn{5}$. In the last case $C$ would lie on a rational normal scroll, or on the Veronese surface.\\
If $C\sub S_4$, where $S_4$ is a scroll of degree $4$ in $\pn{5}$ with hyperplane section $\cH$ and ruling $\cF$, then $C=2\cH+b\cF$ on $S_4$ and thus $S$ is a conic bundle with $K_S\cdot\cH=n-12$. But $(K_S+\cH)^2$ gives $K_S^2=24-3n$, which implies the existence of $3n-16$ singular fibers, i.e. $(-1)$-lines in $S$.\\
Assume now $S_4$ is a Veronese surface, then $S\sub V$,  where $V$ is the cone over $S_4$.  Let $\ol{V}$ be the blowup of $V$ in the vertex, then $\ol{V}$ is a $\pn{1}$-bundle over $\pn{2}$. In this case $S\in|2\cH+b\cF|$ on $\ol{V}$ where $\cF$ is the pullback of a line from $\pn{2}$. Consider $[(\cH-2\cF)^2]\cdot S=-2b$, where $(\cH-2\cF)$ is the contracted divisor.  Since $S$ is smooth $b=0$ and $S$ is a Del Pezzo surface of degree $n=8<11$ in $\pn{6}$ (cf. \ref{nl2}).\\
We can therefore assume that $C$ is not contained in a surface of minimal degree in $\pn{5}$.\\ 
The Hodge index theorem, i.e. $k\leq \frac{e^2}{n}$ implies 
$$n^4-32n^3+332n^2-1120n\leq 8e^2+e(3n^2-80+480)$$
and the bound $e<\frac{n^2}{5}-2n$ yields $ n\leq 27$ and the following list:\\
\begin{center}
\begin{tabular}{|c||c|c|c|c|} \hline
&$n $ & $e$& $k$& $c$\\ \hline\hline

   (a)&  $ 12$&$-2$&$ -3$&$ 3$ \\\hline
   (b)&   $ 12$&$0$&$ -2$&$ 14$\\\hline
  (c)& $ 12$&$ 2$&$-1$&$ 25$\\\hline
  (d)&  $ 12$&$4$&$ 0$&$36$ \\ \hline
  (e)&  $ 14$&$ 0$&$ 0$&$ 0 $ \\\hline
 (f)& $ 16$&$16$&$16$&$80 $ \\ \hline
  (g)& $20$&$ 40$&$70$&$ 206$ \\ \hline
\end{tabular}
\end{center}
Let us examine the various cases:\\
\noindent \underline{case (a)}\\
Since $\chi({\cal O}_S)=0$, $K_S^2=-3$ and $\cH_S\cdot K_S=-2$, the surface $S$ is the blow up of an elliptic ruled surface in $3$ points.\\
Since $S$ has no lines the exceptional curves have degree at least 2.  Now
$(K_S+\cH_S)^2$=$h^0(K_S+\cH_S)=5$, so by Reider's criterium $S$ is embedded in $\pn{4}$ via $|K_S+\cH_S|$. But there are no nonminimal elliptic ruled surfaces of degree 5 in $\pn{4}$, so this case does not occur.\\
\noindent\underline{case (b)}\\
This is a rational surface with $K_S^2=-2$. The adjoint bundle embeds $S$ as a surface of degree $10$ and genus $5$ in $\pn{6}$.  Assuming that $S$ has no lines $|2K_S+\cH_S|$ will blow down $(-1)$- conics and embedd the blown down surface as a Del Pezzo surface of degree 4 in $\pn{4}$.  Therefore $K_S+\cH_S=6l-\sum_1^5 2E_i-\sum_6^{11} E_i$ and $\cH_S=9l-\sum_1^5 3E_i-\sum_6^{11}2 E_i$.  This is example \ref{nl4} in section \ref{noline}.\\
\underline{case (c)}\\
This is the blow up of a $K3$ surface in one point, $\pi:S\to\ol{S}$ and $\cH_S=\pi^{\ast}(\cH_{\ol{S}})-2E$ where $\cH_{\ol{S}}^2=16$.\\
Using the exact sequence:
$$0\to {\cal I}_S\otimes {\cal O}_{\pn{6}}(2)\to {\cal O}_{\pn{6}}(2)\to {\cal O}_S(2)\to 0$$ and the fact that $h^1(S,{\cal O}_S(2))=0$ and thus $h^0(S,{\cal O}_S(2))=24$
we see that $h^0({\cal I}_S\otimes {\cal O}_{\pn{6}}(2))\geq 28-24=4$ and therefore $S$ lies on at least $4$ quadric hypersurfaces.\\
If ${\cal Q}_1\cap ... \cap {\cal Q}_4=S\cup S_4$ where $S_4$ is the residual degree $4$ surface then $S_4$ must be a Veronese surface and $S$ is as in example \ref{nl5}.\\
 Assume now that $\cap \cQ_i$ is not a complete intersection, i.e. $\cap \cQ_i=V$, a threefold in $\pn{6}$. Consider a general $\pn{4}\sub\pn{6}$, then $C=V\cap \pn{4}$ is a curve in $W=\cap\cQ_i\cap \pn{4}$. This only occurs when $C$ lies on a cubic scroll in $\pn{4}$. It is easy to see that this happens only if $V$ lies on a cubic scroll in $\pn{6}$. 
In this case $S$ lies in a cubic $4$-fold scroll in $\pn{6}$. Let $\cF$ be the general fiber and $\cF_S=\cF\cap S$. $\cF_S$ is a pencil of smooth curves on $S$ and $\cF_S^2=0$ or$\cF_S^2=1$. Moreover since $\cF_S\sub\pn{3}$ we know that $\cF$ has no trisecant lines only if the degree$(\cF_S)\leq 4$. Since $S$ is a non-minimal $K3$-surface it cannot have a pencil of rational curves, so $\cF_S$ is an elliptic quartic curve. Since $S$ has degree 12 it must be the complete intersection of the cubic scroll and two quadrics, but this is not $K3$ cf. example \ref{nl5}, so this is not possible.\\
\underline{case (d)}\\
Since $K_S^2=0$ and $\chi(S)=3$ this is an elliptic surface of degree $12$ and $K_S\cdot\cH_S=4$. Let $K_S=mF+\sum(m_i-1)F)$, where $F$ is the general fiber of the elliptic fibration $S\to B$ and $m_iF$ the multiple fibers. 
From $K_S\cdot \cH_S=4$ and  $m=\chi({\cal O}_S)+2g(B)-2\geq 1$ we see that the only possibility is $m=1$ and thus
$g(B)=0$, $K_S=F$ and no multiple fibers. So $|K_S|$ gives a fibration over $\pn{1}$.  The canonical curves are elliptic quartic curves, they each span a $\pn{3}$, and these $\pn{3}$s generate a cubic scroll.  The surface have degree 12 and the canonical pencil have no basepoints so $S$ must be complete intersection of the scroll and two quadric hypersurfaces. This is example \ref{elliptic} in section \ref{noline}. \\
\underline{case e)}\\  Since $\chi(S)=K^2_S=c_2(S)=0$ the surface $S$ must be minimal abelian or bielliptic.  Following Serrano's analysis cf. \cite{Ser} of ample divisors on bielliptic surfaces, one sees that any minimal bielliptic surface of degree 14 in $\pn{6}$ has an elliptic pencil of plane cubic curves, i.e. it has a 3-dimensional family of trisecants. The general abelian surface have no trisecant however, cf \cite{ab1, ab2}.  This is example \ref{nl8}.\\
\underline{case f)}\\
The surface must be the complete intersection of $4$ quadric hypersurfaces in $\pn{6}$, example \ref{nl9} in section \ref{noline}.\\
\underline{case g)}\\
In this case $g(\cH)=p_1$ and thus by (b) in \ref{HA} if the general hyperplane section $C$ is not contained in a minimal surface, it lies on a surface $S_5$ of degree $5$ in ${\bf P}^5$ i.e an anticanonically embedded Del Pezzo surface or the cone over an elliptic quintic curve in $\pn{4}$.\\
In either case the sectional genus of $C$ implies that $C$ is the intersection of $S_5$ with a quartic hypersurface, and each line on $S_5$ will be a 4-secant line to $C$. This excludes case (g).\\

The previous results summarizes as follows:
\begin{theorem} Let $S$ be a smooth surface embedded in $\pn{6}$ with no lines. Then $S$ has no trisecants if and only if it belongs to one of the cases listed in the table below:
\begin{center}
\begin{tabular}{|c|c|c|c|} \hline
surface &degree&linear system& example \\ \hline\hline
${\bf P}^2$&4&${\cal O}_{\bf P}^2(2)$& \ref{nl1} \\\hline
$Bl_7({\bf P}^2)$&8&$6l-\sum 2E_i$&\ref{nl2}\\\hline
$K3$&8&complete intersection (1,2,2,2)&\ref{nl3}\\\hline
$K3$&10&nontrigonal of genus 6&\ref{nl3},\ref{nl7} \\ \hline
$Bl_{11}({\bf P}^2)$&12&$9l-\sum_1^53E_i-\sum_1^6 2E_i$&\ref{nl4} \\ \hline
$Bl_1(K3)$&12&$p^*(\ol{\cH})-2E$&\ref{nl5}\\ \hline
Elliptic&12& minimal with $p_g=2$&\ref{elliptic}\\ \hline
Abelian&14& (1,7) polarization&\ref{nl8}\\ \hline
general type&16& complete intersection (2,2,2,2)&\ref{nl9}\\ \hline
\end{tabular}
\end{center}
\end{theorem}

  %%%%%%%%%%%%%%%%%%%%%%%%%%%%%%%%%%%%%%%%%%%%%%%%%%%%%%%%%%%%%%%%%%%%%%%%%%%%%
%%%%%%%%%%%%%%%%%%%  LIST WITH LINE %%%%%%%%%%%%%%%%%%%%%%%%%%%%%%%%%%%%%%%
%%%%%%%%%%%%%%%%%%%%%%%%%%%%%%%%%%%%%%%%%%%%%%%%%%%%%%%%%%%%%%%%%%%%%%%%%%%%%%%
\section{\Large List of  surfaces with lines}\label{list2}
For surfaces with lines and no proper trisecants the formula (\ref{B}) for the number of tangential trisecants is not necessarily 0.  If $S$ has finitely many lines every $(-1)$-line contributes with multiplicity 1 to the formula.  Therefore we cannot use this formula as it stands in this case. We will use several alternative approaches instead, they will recover all our examples, but will not completely treat all cases, as is made precise in theorem 0.0.1.  One approach we will use is a projection to $\pn{4}$ to get some numerical relations replacing (\ref{B}).
\begin{lemma}\label{lemma} Let $S$ be a smooth surface in $\pn{6}$ with no proper trisecant lines. Assume that $L\subset S$ is a line on the surface, and let $\pi_L:S\to \pn{4}$ be the projection of $S$ from $L$, i.e. the morphism defined by $|\cH_S-L|$.  Then $\pi_L$ is the composition of the contraction of any line on $S$ which meet $L$ and an embedding.  In particular, if $S$ has finitely many lines, $\pi_L(S)$ is smooth unless there are some line $L_1$ on $S$ meeting $L$ with $L_1^2\leq -2$.
\end{lemma}

\proof  Let $P$ be a plane in $\pn{6}$ which contains $L$.  Let $Z_P$ be residual to $L$ in $S\cap P$.  If $Z_P$ is finite, its degree is at most 1, otherwise $S$ would have a trisecant in $P$.  If $Z_P$ contains a curve, this curve would have to be a line which would coincide with $Z_P$, since again $S$ has no trisecants. \qed

With this, let us first examine the cases of surfaces with at least a one dimensional family of lines.
\begin{proposition} Let $S$ be a scroll in $\pn{6}$. Then $S$ has no proper trisecant lines if and only if $S$ is 
\begin{itemize}
\item a rational normal scroll;
\item an elliptic normal scroll, with minimal self intersection $E_0^2=e=1$, cf \ref{ell} .
\end{itemize}
\end{proposition}

\proof If $S$ is a rational normal scroll or elliptic normal scroll with $e=1$, then it is cut out by quadrics and therefore has no trisecant lines.\\
Assume now $S$ is a scroll with no trisecant lines and let $L$ be a line on it. Consider the projection of $S$ from $L$:
$$\pi_L:S\to P^4$$
If there is no line on $S$ meeting $L$ then $\pi_L$ is an embedding by \ref{lemma}.
Since the only smooth scrolls in $\pn{4}$ are the rational cubic scrolls and the elliptic quintic scrolls, $S$ must be as in the statement.\\
If $\pi_L$ is not a finite map, i.e. there is a line section $L_0$, then $S$ must be rational and normal. Assume in fact $S$ not normal, i.e. $S=\pn{}(\Os{}(1)\oplus\Os{}(b))$ with $b\geq 5$, and consider the projection:
$$\pi_{L_0}:S\to \pn{4}$$
The image curve $\pi_{L_0}(S)$ is a rational non normal curve in $\pn{4}$ and therefore it has a trisecant line, $L_t$. Then the linear span $\pn{3}=<L_t,L_0>$ contains three rulings of $S$, therefore also a pencil of  trisecant lines to $S$.\qed
If $S$ is not a scroll, it has only a finite number of lines.
Notice that only lines with self intersection $(-1)$ contributes to the formula (\ref{B}) in \cite{leb}. Therefore we will assume that $S$ has at least one $(-1)$-line, $L$, on it and examine separately the following overlapping cases:
\begin{itemize}
\item $L$ is isolated, i.e. it does not intersect other lines of negative self-intersection.
\item $L$ intersects some $(-1)$-lines.
\item $L$ can be contracted so that it is the exceptional line of an inner projection from $\pn{7}$.
\end{itemize}
%%%%%%%%%%%%%%%%%%%%%%%%%%%%%%%%%%%%%%%%
\subsection{ \large Surfaces with at least one isolated line}
Assume now that $S$ has a finite number of $(-1)$-lines and at least one $L$ is isolated, i.e. it does not intersect any other line. Then the projection:
$$\pi_L:S\to \pn{4}$$ is an embedding by \ref{lemma}. Let $\pi_L^*(\cH)+L=\cH_S$, and notice that $K_{\pi_L(S)}=K_S$.  Using the same invariants $n,e,k,c$ for $S$ we have that
\begin{itemize}
\item degree$(\pi_L(S))=n-3$;
\item $K_{\pi_L(S)}\cdot\cH=e+1$;
\item $K_{\pi_L(S)}^2=k$;
\item $c_2(\pi_L(S))=c$;  
 \end{itemize}
and thus the double point formula in $\pn{4}$ gives
\begin{equation}\label{D}
(n-3)(n-13)-5e-k+c+29=0
\end{equation}
Proceeding as in the previous section, i.e. using the assumptions:
\begin{itemize}
\item (\ref{A}) and (\ref{D});
\item $g\leq p(5)$;
\item $c+k=12\cdot$integer;
\item $k\leq 3c$ and $kn\leq e^2$;
\end{itemize}
numerical computations give the following set of invariants:
\begin{center}
\begin{tabular}{|c||c|c|c|c|} \hline
   &$n $ & $e$& $k$& $c$\\ \hline\hline
(a)&  $ 8$&$-8$&$ 5$&$ -5$ \\\hline
   (b)&   $ 8$&$-4$&$ 1$&$ 11$\\\hline
  (c)&  $ 9$&$ -3$&$ -1$&$ 13 $ \\\hline
(d)&  $ 10$&$ -2$&$ -2$&$ 14 $ \\\hline
(e)&  $ 11$&$ 1$&$ -1$&$ 25 $ \\\hline
\end{tabular}
\end{center}
\noindent \underline{case (a)}\\
$\chi(\Os{S})=0$ and $K_S^2=5$, this is impossible.\\
\noindent\underline{cases (b), (c), (d)}\\
Since $e=\cH_S\cdot K_S\leq -2$ and $\chi(\Os{S})=1$ the surfaces $S$ are rational.  The adjunction morphism defined by $|\cH_S+K_S|$ is birational and maps $S$ onto $\pn{2}$, a quadric in $\pn{3}$ and a Del Pezzo surface of degree 4 in $\pn{4}$ respectively in the three cases. Thus we recover the surfaces constructed in \ref{p48}, \ref{p49} and \ref{p410}.\\
\noindent\underline{case(e)}\\
Since $\cH_S\cdot K_S=1$ and $\chi(\Os{S})=2$ the surface $S$ must be a $K3$-surface blown up in one point. Thus we recover the inner projection of the general $K3$-surface of degree $12$ in $\pn{7}$ described in \ref{prk3}. The line $L$ is the only exceptional line on the surface.\\
We have then proved:
\begin{proposition}\label{proj4}
Let $S$ be a surface in $\pn{6}$ with no trisecant lines and with at least one isolated $(-1)$-line. then $S$ is one of the following:
\begin{itemize}
\item A rational surface of degree $8$ and genus $3$, as in \ref{p48};
  \item A rational surface of degree $9$ and genus $4$, as in \ref{p49};
\item A rational surface of degree $10$ and genus $5$, as in \ref{p410};
\item A non minimal $K3$-surface of degree $11$ and genus $ 7$, as in \ref{prk3}.
\end{itemize}
\end{proposition}
%%%%%%%%%%%%%%%%%%%%%%%%%%%%%%%%%%%%%%%%%%%%%
\subsection{ \large Conic bundles}\label{conic}
%%%%%%%%%%%%%%%%%%%%%%%%%%%%%%%%%%%%%%%%%%%%%%%%%%%%%%%%%%%%%%%%%%%%%%%%%
Assume now that $S$ has at least two $(-1)$-lines $L_1$ and $L_2$ which meet.\\  Then $(L_1+L_2)^2=0$ and $L_1+L_2$ or some multiple of it moves in an algebraic pencil of conic sections. Thus $S$ is a conic bundle.\\ 
If there is a line $L$ of selfintersection $L^2\leq -2$ intersecting $L_1$, then this line intersects all the members of the pencil and therefore is mapped onto the base curve. This would imply then that $S$ is a rational conic bundle in $\pn{6}$. Moreover the general hyperplane section $C\in|\cH_S|$ is a hyperelliptic curve which, since we are in the smooth case, must be nonspecial. This implies $h^1(S,\cH_S)=0$ and thus $7=1+\frac{n-e}{2}$ by Riemann Roch. Since $(K_S+\cH_S)^2=n+2e+k=0$ we derive:
\begin{itemize}
\item $e=n-12$;
\item $k=24-3n$;
\item $c=3n-12$.
\end{itemize}
Substituting those values in the first Le Barz equation, (\ref{A}), we get an equation of degree $3$ in the variable $n$ and therefore three possible cases, which are easily seen to be the ones constructed in \ref{con1}, \ref{con2} and \ref{con3}.\\
Let us now assume that there is no line $L$ of selfintersection $L^2\leq -2$ intersecting $L_1$. Then the projection:
$$\pi_{L_1}:S\to \pn{4}$$ will be composed by the contraction of $L_2$ and possibly other $-1$-lines intersecting $L_1$, and an embedding by \ref{lemma}.  The conic bundle structure is of course preserved. Therefore the image surface in $\pn{4}$ is rational of degree 4 or 5, or it is an elliptic conic bundle of degree $8$ cf. \cite{BR, ES}. It is clear that the two first come from the surfaces of type \ref{con2} and \ref{con3}.  The elliptic conic bundle has two plane quartic curves on it in $\pn{4}$ which are bisections on the surface.  On $S$ their preimage would have degree $5$ or $6$ depending on the intersection with $L_1$.  But this is a curve of genus $3$ so the degree upstairs must be $6$ if $S$ is smooth.  A curve of degree $6$ and genus $3$ has trisecants, so this excludes this case.
%%%%%%%%%%%%%%%%%%%%%%%%%%%%%%%%%%%%%%%%%%%%%%%%%%%%%%%%%%%%%%%%%%%%%%%%%
\subsection{ \large Inner projections from $\pn{7}$}
%%%%%%%%%%%%%%%%%%%%%%%%%%%%%%%%%%%%%%%%%%%%%%%%%%%%%%%%%%%%%%%%%%%
In this section we will investigate the remaining cases, i.e. surfaces with $r$
 $(-1)$-lines on it, which do not intersect.
In other words $S$ has $r$ exceptional curves of the first kind $E_1,...,E_r$, which can
 therefore be contracted.\\
 Let us assume that we can contract at least one of the $E_i$'s down
  to  $x\in\ol{S}$ which is not a base point for $\ol{\cH}$, where $\cH=p^*(
  \ol{\cH})-E$ and $p$ is the projection map $p :S\to \ol{S}$. Then we can
  assume that $\ol{\cH}$ embeds $\ol{S}$ in $\pn{7}$ with no trisecant lines
  and with $(r-1)$ $(-1)$-lines.\\
  Notice that the existence of $r$ $(-1)$-lines on $S$ gives a contribution of $4r$ in
  the number of tangential trisecants, which in terms of the formula as stated
  in (\ref{B}) means: 
  $$ T_3=4r$$
  Moreover the existence of $(r-1)$ $(-1)$-lines on $\ol{S}$ gives a
  contribution of $-(r-1)$ in the formula for trisecant lines of surfaces in
  $\pn{7}$. Using the same invariants $n,e,k,c$ for $S$, as in \ref{nlist},  we have:
\begin{itemize}
\item degree($\ol{S}$)=$n+1$;
\item $K_{\ol{S}}\cdot\ol{\cH}=e-1$;
\item $K_{\ol{S}}^2=k+1$;
\item $c_2(\ol{S})=c-1$;
\end{itemize}
Plugging those invariants in the formula of Le Barz we get:
\begin{equation}\label{C}
S_3:=n^3-27n^2+176n+108+c(3n-37)-k(3n-53)-e(15n-177)=6r-6
\end{equation}
We can then assume:
\[ \left\{\begin{array}{ll}
D_3&=2 n^3  - 42 n^2  + 196 n-K(3n-28)+c(3n-20) - e(18n - 132)=0\\
T_3&=6 n^2  - 84 n+k(n-28) -c(n-20)+ e(4n-84)=4r\\
S_3&=n^3-27n^2+176n+108+c(3n-37)-k(3n-53)-e(15n-177)
\\&=-6r+6
\end{array} \right. \]
Thus \[ \begin{array}{ll}
2S_3&=-2(6r-6)\\
&=D_3-12n^2+156n+108-k(3n-78)+c(3n-54)-e(12n-222)\\
&=D_3-3T_3+6n^2-96n+216-6k+6c-30e\\
&=-12r+6n^2-96n+216-6k+6c-30e
\end{array}  \]
and hence 
\begin{equation} \label{ss}12=6n^2-96n+216-6k+6c-30e\end{equation}
Let $\ol{C}$ be a general hyperplane 
section of $\ol{S}$, of degree $n+1$ and genus $g$. If $n+1\geq 15$, by \ref{p3}
 either $g\leq p_2$ or  $\ol{C}$ lies on a surface of degree $5$ or $6$ in $\pn{6}$. If $\ol{C}\sub S_5$, where $S_5$ is a surface of degree $5$ in $\pn{6}$, then $S_5$ is a rational normal scroll and we can write $\ol{C}=2\cH+b\cF$, since we are assuming $\ol{C}$ has no trisecant lines. Like in section \ref{nlist}, $S_5$ can be viewed as a conic bundle  and thus this case has already been studied in \ref{conic}.\\
If $\ol{C}\sub S_6$, where $S_6$ is a surface of degree $6$ in $\pn{6}$,
 then $S_6$ is a Del Pezzo surface embedded via the anticanonical bundle 
 $-K_{S_6}$, i.e. $S_6=Bl_3(\pn{2})$, or it is the cone over an elliptic curve of degree 6 in $\pn{5}$. On the cone the curve $\ol{C}$ would intersect any ruling at most twice.  Since it is smooth it passes through the vertex at most once, so the degree $n+1\leq 13$.  In case of a Del Pezzo surface, write $\ol{C}=al-a_1E_1-a_2E_2-a_3E_3$, 
then the fact that $\ol{C}$ has no trisecants yields:
$a_i\leq 2$, $i=1,2,3$ and $a-a_i-a_j\leq 2$, which implies $a\leq 6$. Since 
$n\geq 14$ the only possibility would be $a=6$ and $a_i<2$ for some $i$, but this is impossible since we
need $a-a_i-a_j\leq 2$. \\
Thus we can assume $g\leq p_3$ (cf. \ref{p3}) when $n\geq 13$.
Then easy computations involving (\ref{ss}), (\ref{A}), $kn\leq e^2$ (Hodge Index theorem) yield $n\leq 15$.\\
This bound and
\begin{itemize}
\item(\ref{A}), (\ref{B}), (\ref{C});
\item $k\cdot n\leq e^2$ and $k\leq 3c$;
\item $n+e+2\leq p(5)$ (Castelnuovo bound);
\item $n+2e+k=(K_S+\cH_S)^2>0$ (not a conic bundle case).
\end{itemize}
give the following numerical invariants:
\begin{center}
\begin{tabular}{|c||c|c|c|c|c|} \hline
&$n $ & $e$& $k$& $c$&r\\ \hline\hline
(a)&  $ 8$&$-4$&$ 1$&$ 11$&$8$ \\\hline
(b)&   $ 9$&$-3$&$ -1$&$ 13$&$9$\\\hline
(c)&  $ 10$&$-2$&$ -2$&$14$&$6$ \\ \hline
(d)&  $ 11$&$ 1$&$ -1$&$ 25 $ &$1$\\\hline
\end{tabular}
\end{center}
One can immediately see that we recovered the examples \ref{p48}, \ref{p48}, \ref{p410} and \ref{prk3} respectively.

The results of this section summarize as follows:
\begin{proposition}
Let $S$ be a surface embedded in $\pn{6}$ with  no trisecant lines and having 
$r$ $(-1)$-lines on it. Assume that $S$ is the inner projection of a smooth surface $\ol{S}\sub\pn{7}$, then $S$ is as in Proposition \ref{proj4} i.e. 
\begin{itemize}
\item a rational surface of degree $8$ and genus $3$ with eight $(-1)$-lines, as in \ref{p48};
  \item a rational surface of degree $9$ and genus $4$ with nine $(-1)$-lines, as in \ref{p49};
\item a rational surface of degree $10$ and genus $5$ with six $(-1)$-lines, as in \ref{p410};
\item A non minimal $K3$-surface of degree $11$ and genus $7$ with only one $(-1)$-line, as in \ref{prk3}.
\end{itemize}
\end{proposition}
%%%%%%%%%%%%%%%%%%%%%%%%%%%%%%%%%%%%%%%%%%%%%%%%%%%%%%%%%%%%%%%
%%%%%%%%%%%%%%%%%%%%%%%%%%%%%%%%%%%%%%%%%%%%%%%%%%%%%%%%%%%%%%%%%
\section*{\vskip-1cm \Large Conclusion}
We may summarize sections  \ref{nlist}, \ref{list2} in the following
\begin{proposition}
 Let $S$ be a surface in $\pn{6}$ with no trisecant lines.
Unless $S$ is not an inner projection from $\pn{7}$ and every $(-1)$-line on $S$, if there are any, meet some other line $L$ on the surface with $L^2\leq -2$, the surface belongs to the list of examples in sections \ref{noline} and \ref{wline}. 
\end{proposition}
We have made computations as in section \ref{nlist} fixing the number $r$ of $(-1)$-lines on the surface. Checking up to $r=100$ give no contribution to the list we have produced.\\ 
This numerical observation and the fact that the examples constructed in section \ref{wline} cover all the cases listed in section \ref{list2} lead us to make  the following
\begin{con}
Let $S$ be a surface in $\pn{6}$ with no trisecant line, then
the surface belongs to the list of examples in sections \ref{noline} and \ref{wline}. 
\end{con}

%%%%%%%%%%%%%%%%%%%%%%%%%%%%%%%%%%%%%%%%%%%%%%%%%%%%%%%%%%%%%%%%%%%%%%%%%%%%%%%%%%%%%%%%%%%%%%%%

%%%%%%%%%%%%%%%%%%%%%%%%%%%%%%%%%%%%%%%%%%%%%%%%%%%%%%%%%%%%%%%%%%%%%%%%%%%
\bigskip
 Sandra Di Rocco\\
{\it KTH, the Royal Institute of Technology \\        
Department of Mathematics  \\          
S-100 44 Stockholm, SWEDEN\\}
e-mail: sandra@math.kth.se\\ \\
Kristian Ranestad\\
{\it Matematisk Institutt, UiO\\
 P.B. 1053 Blindern,\\
N-0316 Oslo, NORWAY}\\
e-mail: ranestad@math.uio.no\\

\end{document}